
\input amstex
\documentstyle{amsppt}
\pagewidth{125mm}
\pageheight{195mm}

\topmatter
\title
On one family of 13-dimensional closed Riemannian
positively curved manifolds
\endtitle
\author
Ya. V. Bazaikin
\endauthor
\affil
Institute of Mathematics, 630090 Novosibirsk, Russia
\endaffil
\endtopmatter
\leftheadtext{}
\rightheadtext{}
\document
\NoBlackBoxes

\head
1. Introduction and main results
\endhead

'In the present paper we describe one family of closed Riemannian manifolds
with
positive sectional curvature.

ŠNow the list of known examples is not large (for instance,
all known manifolds with
dimension $> 24$ are diffeomorphic to compact rank one symmetric spaces)
(we restrict ourselves only by pointing out simply connected manifolds) :

1) Berger described all normally homogeneous closed positively curved
manifolds that are compact rank one symmetric spaces (i.e., the spheres $S^n$,
the complex projective spaces $CP^n$, the quaternionic projective
spaces $HP^n$, and the projective Cayley plane $CaP^2$), and two exceptional
spaces of form $Sp(2)/SU(2)$ and $SU(5)/Sp(2) \times S^1$ with dimension
7 and 13,  respectively (notice that the embedding $SU(2) \subset Sp(2)$ is not
standard) (\cite{Be});

2) "Wallach had shown that all even-dimensional simply connected closed
positively curved are  diffeomorphic to normally
homogeneous ones or the flag spaces over $CP^2, HP^2$, and $CaP^2$
(with dimension 6,12, and 24, respectively)
(\cite{W});

3) €Aloff and Wallach (\cite{AW}) constructed infinite series of spaces
$N_{p,q}$ of the form  $SU(3)/S^1$ where the subgroup $S^1$ is a winding of a
maximal torus of group $SU(3)$ and, since that, is defined by a pair
of relatively prime integer parameters $p$ and $q$. If some conditions for
these parameters $p$ and $q$ hold then these manifolds admit left-invariant
homogeneous Riemannian metric with positive sectional curvature.
Berard-Bergery (\cite{BB}) had shown that the Aloff-Wallach spaces are all
possible manifolds that admit homogeneous positively curved metric and do not
admit normally homogeneous one, and Kreck and Stolz had found
among them a pair of homeomorphic but nondiffeomorphic manifolds
($N_{-56788,5227}$ and $N_{-42652,61213}$ ; \cite{KS}) ;

4) by using of the construction of Aloff and Wallach,
Eschenburg had found an infinite
series of seven-dimensional spaces with nonhomogeneous
positively curved metrics (\cite{E1}) and in the sequel had found an example of
six-dimensional nonhomogeneous space with positively curved metric
(\cite{E2}).

This list contains all known up to now topological types
of simply connected closed manifolds that admit metrics with positive
sectional curvature. Notice that only two of them have dimension 13 :
the sphere $S^{13}$ and the normally homogeneous Berger space
$SU(5)/Sp(2)\times S^1$.

ŽThe main result of the present paper is the construction of the new series of
simply connected closed 13-dimensional manifolds that admit positively curved
metrics.
In particular, we prove the following theorem.

\proclaim{Main Theorem}
 Let $U(5)$ be a group of complex unitary $5\times 5$-matrices and a group
 $U(4)\times U(1)$ is embedded into it as a subgroup of matrices of  block form
with two blocks with size $4\times 4$ and $1\times 1$. Let $M^{25}$ be a
homogeneous
Riemannian manifold diffeomorphic to  $U(5)$ and endowed by metric induced
from two-sided invariant metric on
 $U(5) \times U(4)\times U(1)$ by projection
$$
U(5) \times U(4) \times U(1) \rightarrow U(5) \times U(4) \times U(1) /
U(4) \times U(1) = M^{25}
$$
with diagonal embedding $U(4) \times U(1) \rightarrow U(5) \times U(4) \times
U$ ($g \rightarrow (g,g) \in U(5) \times (U(4) \times
U(1))$).

…Let $\bar p = (p_1,\dots,p_5)$ be a 5-tupel of integer numbers
such that for every
permutation
 $\sigma \in S_5$
the following conditions hold

a)  $p_{\sigma (1)}+p_{\sigma (2)}-p_{\sigma (3)}-p_{\sigma (4)}$
¢§ ¨¬­® ¯à®áâ® á $p_{\sigma (5)}$ ,

b) $p_{\sigma (1)}+p_{\sigma (2)}+p_{\sigma (3)}>p_{\sigma
(4)}+p_{\sigma (5)}$ ,

c) $p_{\sigma (1)}+p_{\sigma (2)}+p_{\sigma (3)}+p_{\sigma
(4)}>3p_{\sigma (5)}$ ,

d) $3(p_{\sigma (1)}+p_{\sigma (2)})>p_{\sigma (3)}+p_{\sigma
(4)}+p_{\sigma (5)}$ .

Let $M_{\bar p}$ be a  factor-space $M_{\bar p}$ of
$M^{25}$ under the action of
$S^1 \times (Sp(2) \times S^1)$ given by
$$
(z_1,(A,z_2)): X \rightarrow
diag(z_1^{p_1},z_1^{p_2},z_1^{p_3},z_1^{p_4},z_1^{p_5})\cdot
X\cdot \left( \matrix A^{*} \bar{z}_2 & 0 \\  0 & 1
\endmatrix \right),
$$
where $X\in M^{25}$, $z_1, z_2\in S^1$, and $A\in Sp(2)$. Let
$M_{\bar p}$ is endowed by metric induced by factorization
$M^{25} \rightarrow M_{\bar p}$ then

1)  $M_{\bar p}$ is simply connected and  $\dim M_{\bar p} = 13$ ;

2)  $M_{\bar p}$ has positive sectional curvature ;

3) the groups of cohomologies of  $M_{\bar{p}}$ are the following ones :

$$
H^i=
\cases
 {\bold Z}, & for\ \ \  $i=0,2,4,9,11,13,$ \\
 0, & for\ \ \  $i=1,3,5,7,10,12;$
\endcases
$$
the groups  $H^6$ and $H^8$ are finite and their orders are equal to
$|\sigma_1^3-4\sigma_1\sigma_2+8\sigma_3|$ where $\sigma_k$ is a value of an
elementary symmetric polynom, of degree $k$, on five
variables at  $(p_1,\dots,p_5)$.
\endproclaim

{\bf ‡ Remarks.}
"The condition a)-d) hold, for instance, for $p_1 =1,
p_2=p_3=p_4=p_5=q^n$ where  $q$ is a prime number. In this case
the order of $H^6(M_{\bar p})$ is equal to $r(q,n) = 8q^{2n}-4q^{n}+1$ and
$r(q,n) \rightarrow \infty$ as $q \rightarrow \infty$ . It follows from
this example that there exists infinitely many pairwise nonhomeomorphic
closed simply connected positively curved manifolds of the form $M_{\bar p}$.

ŒOne can see that for $n=0$ we obtain manifold which is diffeomorphic to the
13-dimensional Berger space.

'There exist another series and the simplest construction of them
was pointed out to us by U. Abresch. In particular, let take
a 5-tupel of numbers for which the condition a) of Main Theorem holds (notice
here that we call two numbers relatively prime if their maximal common
divisor is equal to one) and no one of numbers
$|p_{\sigma(1)}+p_{\sigma(2)}-p_{\sigma(3)}-p_{\sigma(4)}|$
vanishes. Let add to  $p_1,\dots,p_5$ the same natural number
$a_n = n \cdot \prod_{\sigma \in S_5}
|p_{\sigma(1)}+p_{\sigma(2)}-p_{\sigma(3)}-p_{\sigma(4)}|$.
‹One can see that there exists sufficiently large number $N$ such that
for all $n > N$ 5-tupels
$(p_1+a_n,\dots,p_5+a_n)$ satisfy to conditions b)-d) and it is easy to observe
that these tupels always satisfy condition a). For instance,
one can start from initial 5-tupel of the form
$(1,1,1,2q,4q)$ where
$q$ is a prime number.

General structure of the family of constructed (in Theorem 1) manifolds and
pinchings of their metrics will be considered separately.

At the construction of metric we follow to methods developed in \cite{E1}.
But at the proving of positivity of curvature the method introduced in
\cite{E1}
met some difficulties which are passed by using of Lemma 8.

'In the next chapter the spaces $M_{\bar p}$ are constructed, the statement 2
of
the theorem is proved in chapter 3 (Theorem 1), and statements 1 and 3
are proved in the  fourth, final, chapter (Theorems 2 and 3, respectively).

 This work is supported by the Russian Foundation for Fundamental
Researches (grant No 94-01-00528).

€Author thank his adviser I.A. Taimanov for posing the problem and helpful
advices.  Author thank U. Abresch for helpful comments.

\head
2. Construction of spaces $M_{\bar p}$
\endhead

\subhead
2.1. Riemannian submersion and it's properties
\endsubhead

Let $M$ and $N$ be Riemannian manifolds and  $f: M\rightarrow N$ be a smooth
mapping. A mapping $f$ is called {\it submersion} if $f$ is surjective
(i.e., $f(M)=N$) and the linear mapping $d_x:T_xM\rightarrow
T_{f(x)}N$ is isomorphism for every point $x \in M$. Then at every point
$x \in M$ the tangent space to $M$ canonically decomposes into a direct sum of
two subspaces $T_xM=(T_xM)^v\oplus (T_xM)^h$ where
$$
(T_xM)^v=T_xK,\ \
K=f^{-1}(f(x))
$$
 and $(T_xM)^h$ is an orthogonal complement to $(T_xM)^v$.
These subspaces are called vertical and horizontal, respectively. 
It is evident that  $d_xf|_{(T_xM)^h}:(T_xM)^h\rightarrow
T_{f(x)}N$ is an isomorphism.
If this isomorphism preserves metric the mapping $f$ called {\it Riemannian
submersion}.

The next lemma gives a general construction of examples of Riemannian
submersions.

\proclaim{Lemma 1}
Let  $G$ be a group of isometries that acts freely and with closed orbits on a
Riemannian manifold $M$. Then on the space of orbits one can introduce
the structure of Riemannian manifold such that the natural projection
 $\pi: M\rightarrow N$ be a Riemannian submersion.
\endproclaim

'For Riemannian submersions the curvatures of manifolds $M$ and $N$ are related
by
formula found in \cite{ON}. We restrict ourselves only by it's
corrolary that we will need in the sequel.

\proclaim{Lemma 2}
Let $\pi:M\rightarrow N$ be a Riemannian submersion.
Put  $x\in M, y\in N,\pi (x)=y$. If …$\sigma^*$ is a two-dimensional horizontal
plane in $T_xM$ and  $\sigma=d_x\pi(\sigma^*)$ then
$$
K(\sigma)\geq K(\sigma^*).
$$
\endproclaim

"The proof of the next lemma one can find, for instance, in \cite{M}.

\proclaim{Lemma 3}
Let $G$ be a Lie group with two-sided invariant metric
$\langle \ ,\ \rangle$ and ${\bold g}$ be a tangent space, at the unit, endowed
by a structure of Lie algebra. Then for any $X,Y \in {\bold g}$
a sectional curvature in direction $Span(X,Y)$ is equal to
$$
K(X,Y)=\frac{1}{4} \langle [X,Y],[X,Y] \rangle
$$
\endproclaim

\subhead
2.2. Normally homogeneous metric on $U(5)$
\endsubhead

'In this subchapter we construct one Riemannian metric on
$U(5)$ and in the next subchapter we define free actions of the group
$S^1 \times (Sp(2) \times S^1)/Z_2$ on  $U(5)$  that are isometries with
respect to this metric. This construction of an  auxiliary metric on
$U(5)$ gives itself an example of Riemannian submersion. Moreover  metrics
looked for on the spaces of orbits, of action
$S^1 \times (Sp(2) \times S^1)/Z_2$ on­ $U(5)$, will be constructed by this
metric with using of Lemma 1. At this construction we follow paper \cite{E2}.

Let $G$ be the  Lie group  $U(5)$ and  $K=U(4)\times U(1)$ be the subgroup
which is
embedded in the standard manner. Let consider a usual two-sided invariant
Riemannian metric
$\langle\ ,\ \rangle_0$ on $G$ :
$$
\langle X,Y \rangle_0 = Re\ \ trace\ \ (XY^*)
$$
where $X,Y \in {\bold u}(5)$.
This metric canonically induces metrics on
$K$ and $G\times K$. These metrics we will also define
by $\langle\ ,\ \rangle_0$.

Let  $\triangle K=\{ (k,k) | k \in K \}$ be a subgroup in
$G\times K$.
 We consider an action $\triangle K$ on $G\times K$ by right shifts :
$$
((g,k),k') \longmapsto (gk',kk')
$$
for $ g\in G, k,k'\in K $.
ŽIt is evident that that is an isometrical free action.
By Lemma 1, there exists a metric on the space of orbits
$(G\times K)/\triangle K$
such that the natural projection
$$
\pi:G\times K \rightarrow (G\times K)/\triangle K
$$
is a Riemannian submersion.
One can see that the correspondence
$(g,k)\rightarrow gk^{-1}$ gives a diffeomorphism
$(G\times K)/\triangle K$ with $G$. By pulling,
with the help of this diffeomorphism, the Riemannian metric from
the space of orbits
$(G\times K)/\triangle K$ onto
$G$ we obtain a metric $\langle\ ,\ \rangle$ on $G$. 
For that the mapping
$$
\pi:G\times K \rightarrow G:(g,k)\longmapsto gk^{-1}
$$
is a Riemannian submersion.

 Let consider a left shift by element $(g,k^{-1})$ on the group
$G\times K$ where $g\in G, k\in K$. ' Since the metric $\langle\ ,\ \rangle_0$
is two-sided invariant, this mapping is an isometry.
Moreover, the left shift maps fibers of submersion into fibers and
hence induces a mapping
$$
g'\longmapsto gg'k:G\rightarrow G,
$$
on $G$, that is an isometry. Thus,
we conclude that {\it the metric $\langle \ ,\ \rangle$
is left-invariant under $G$  and right-invariant under $K$.}

Let ${\bold k}={\bold u }(4)\oplus {\bold u}(1)$ and ${\bold g}={\bold u}(5)$
be tangent algebras of groups $G$ and $K$ , respectively.
We denote by  Ž${\bold p}$
an orthogonal complement to ${\bold k}$ in ${\bold g}$ with respect to the
metric
$\langle\ ,\ \rangle_{0}$. Then the decomposition ' ${\bold g}={\bold k}\oplus
{\bold p}$ is invariant under  $Ad(K)$. Moreover,  Š $G/K$ is the symmetric
space
$CP^4$ and, since that, we have
$$
[{\bold k},{\bold k}]\subset {\bold k},\ \
[{\bold p},{\bold p}]\subset {\bold k},\ \
[{\bold k},{\bold p}]\subset {\bold p}.
\eqno{(1)}
$$

'The vertical subspace of submersion $\pi$
at $(e,e)$ is
$$
V=\{(Z,Z)|Z\in {\bold k}\}=\triangle {\bold k}.
$$
Hence  $(X,Y)\in {\bold g}\oplus {\bold k}$ lies in the horizontal subspace
$H$ if
$$
\langle(X,Y),(Z,Z)\rangle_0=0
$$
for every $ Z\in {\bold k}$
that implies
$$ \langle X,Z \rangle_0+ \langle Y,Z
\rangle_0=0,
$$
$$
\langle X+Y,Z \rangle_0=0
$$
for every $Z\in {\bold k}$.
‡ Notice that in this case $X+Y\in {\bold p}$, i.e., $X_k+Y_k=0$ and
$Y=Y_k=-X_k$. We derive from that that
$$
H=\{(X_k+X_p,-X_k)|X_k\in {\bold k},\ X_p\in {\bold p}\}
$$
and $d_{(e,e)}\pi |_H:H\rightarrow {\bold  g}$ is an isometry.

' Since $d_{(e,e)}\pi(X,Y)=X-Y$, for every $X\in {\bold g}$
the following equality holds
$$
(d_{(e,e)}\pi |_H)^{-1}(X)=(\frac{1}{2}X_k+X_p,-\frac{1}{2}X_k).
\eqno{(2)}
$$
holds.

\proclaim{Lemma 4}
Let $X\in {\bold g}$ and $Y\in {\bold k}$.  Then $\langle
X,Y\rangle=\frac{1}{2}\langle X,Y\rangle_0$.
\endproclaim

{\bf "Proof of Lemma 4.}

By '(2), we have
$$
\langle X,Y\rangle=
\langle(\frac{1}{2}X_k+X_p,-\frac{1}{2}X_k),
(\frac{1}{2}Y_k+Y_p,-\frac{1}{2}Y_k)\rangle_0=
$$
$$
=\langle \frac{1}{2}X_k+X_p,\frac{1}{2}Y\rangle_0+\langle
-\frac{1}{2}X_k,-\frac{1}{2}Y\rangle_0=
\langle X_k+X_p,\frac{1}{2}Y\rangle_0=\frac{1}{2}\langle
X,Y\rangle_0.
$$
Lemma 4 is proved.

'In the sequel we will mean by curvature of the space $G$ it's
curvature with respect to the metric $\langle\ ,\ \rangle$.

\proclaim{Lemma 5}
Let $\sigma$ be a two-dimensional plane in  ${\bold g}$ and  $K(\sigma)=0$.
'Then $\sigma=Span(X,Y)$, where $X\in {\bold g},\ Y\in {\bold k}$ and
$$
[X_p,Y]=[X_k,Y]=0
$$
\endproclaim

{\bf "Proof of Lemma  5.}

Let $\sigma=Span(X,Y)$ where $X,Y\in {\bold g}$. Let
$\sigma^*=Span((\frac{1}{2}X_k+X_p,-\frac{1}{2}X_p),
(\frac{1}{2}Y_k+Y_p,-\frac{1}{2}Y_k))$ lies in the horizontal
subspace of submersion $\pi$.
ˆWe have $d_e\pi(\sigma^*)=\sigma$. By Lemmas 2 and 3, $0\leq
K(\sigma^*)\leq K(\sigma)=0$. ŽThat implies  $K(\sigma^*)=0$.

ˆIt follows from lemma 3 that
$$
[(\frac{1}{2}X_k+X_p,-\frac{1}{2}X_k),
(\frac{1}{2}Y_k+Y_p,-\frac{1}{2}Y_k)]=0,
$$
$$
([\frac{1}{2}X_k+X_p,\frac{1}{2}Y_k+Y_p],
[-\frac{1}{2}X_k,-\frac{1}{2}Y_k])=0.
$$
' Hence ,
$$
[X_k,Y_k]=0,
$$
$$
\frac{1}{2}[X_k,Y_p]+\frac{1}{2}[X_p,Y_k]+[X_p,Y_p]=0
$$
'By (1),
$[X_p,Y_p]\in {\bold k}$ and $[X_k,Y_p]+[X_p,Y_k]\in {\bold p}$, that implies
$$
[X_p,Y_p]=0,
$$
$$
[X_k,Y_p]+[X_p,Y_k]=0.
$$

" Next,  $X_p,Y_p\in {\bold p}$ are tangent vectors to  positively curved space
$CP^4$. Since the curvature of  $CP^4$ in the direction
$Span(X_p,Y_p)$ vanishes, vectors  $X_p,Y_p$ are linearly dependent.
Hence, we may assume that
$\sigma=Span(X,Y)$ where
$Y\in {\bold k}$.

'Then we obtain
$$
[X_k,Y]=[X_p,Y]=0.
$$

Lemma 5 is proved.

\subhead
2.3. Free actions on $U(5)$ and construction of spaces $M_{\bar p}$
\endsubhead

Let $p_1,p_2,p_3,p_4,p_5$ be integer numbers. Put
$P'=S^1\times(Sp(2)\times S^1)$ where we assume that $Sp(2)$ is standard
embedded into  $SU(4)$.

 Let consider an action of the group $P'$ on $G=U(5)$:
$$
(z_1,(A,z_2)):X\mapsto
diag(z_1^{p_1},z_1^{p_2},z_1^{p_3},z_1^{p_4},z_1^{p_5})\cdot
X\cdot \left( \matrix A^{*} \bar{z}_2 & 0 \\  0 & 1
\endmatrix \right),
$$
where $X\in G$, $z_1,z_2\in S^1$, $A\in Sp(2)$.

\proclaim{Lemma 6}
Let $p_{\sigma (1)}+p_{\sigma (2)}-p_{\sigma (3)}-p_{\sigma
(4)}$ is relatively prime with $p_{\sigma (5)}$ for every transposition
$\sigma \in S_5$. 'Then this action has a kernel isomorphic to
$Z_2=(1,\pm(E,1))$ and therefore induces a free action of group
$$
P=S^1\times \frac{Sp(2)\times
S^1}{\pm (E,1)}=:P_1\times P_2.
$$
on $G$.
\endproclaim

{\bf "Proof of Lemma 6.}

"Let assume that
$$
X=diag(z_1^{p_1},z_1^{p_2},z_1^{p_3},z_1^{p_4},z_1^{p_5})\cdot X
\cdot \left( \matrix A^{*} \bar{z}_2 & 0 \\  0 & 1
\endmatrix \right),
$$
$$
diag(\bar{z}_1^{p_1},\bar{z}_1^{p_2},\bar{z}_1^{p_3},
\bar{z}_1^{p_4},\bar{z}_1^{p_5})=X
\left( \matrix A^{*} \bar{z}_2 & 0 \\  0 & 1
\endmatrix \right) X^{-1}.
$$
'We consider the maximal torus in  $Sp(2)$ :
$$
T^2=\{diag(u,v,\bar{u},\bar{v})|u,v\in S^1\}.
$$
'Then there exists an element  $Y\in Sp(2)$ such that
$A^{*}=Y diag(u,v,\bar{u},\bar{v}) Y^{-1}$ for some  $u,v\in S^1$.
ˆSo,  we have
$$
diag(\bar{z}_1^{p_1},\bar{z}_1^{p_2},\bar{z}_1^{p_3},
\bar{z}_1^{p_4},\bar{z}_1^{p_5})=
$$
$$
=\left( X \left( \matrix Y & 0 \\  0 & 1 \endmatrix
\right) \right) diag(u\bar{z}_2,v\bar{z}_2,\bar{u} \bar{z}_2,\bar{v}
\bar{z}_2,1) \left( X \left( \matrix Y & 0 \\  0 & 1
\endmatrix \right) \right) ^{-1}.
$$
‡That means that there exist a permutation $i$ such that
$\{i_1,i_2,i_3,i_4,i_5\}=\{1,2,3,4,5\}$ and
$$
\bar{z}_1^{p_{i_1}}=u\bar{z}_2,\ \
\bar{z}_1^{p_{i_2}}=v\bar{z}_2,\ \
\bar{z}_1^{p_{i_3}}=\bar{u} \bar{z}_2,\ \
\bar{z}_1^{p_{i_4}}=\bar{v} \bar{z}_2,\ \
\bar{z}_1^{p_{i_5}}=1.
$$
ˆIt follow from the first four equalities that
$$
\bar{z}_1^{p_{i_1}+p_{i_3}-p_{i_2}-p_{i_4}}=1.
$$
By the condition of Lemma,  $p_{i_5}$  is relatively prime with
$p_{i_1}+p_{i_3}-p_{i_2}-p_{i_4}$, and therefore  $\bar{z}_1=1$, i.e.,
$z_1=1$.
'Then  $\bar{z}_2^2=1$, $z_2=\pm 1$.

1) …If $z_2=1$ then  $u=v=1, A^{*}=YEY^{-1}=E$, i.e., $ A=E$.

2) …If $z_2=-1$ then $u=v=-1, A=-E$.

ˆThus, the kernel of action is given by
$$
Z_2=(1,\pm(E,1)).
$$
Lemma 6 is proved.

ƒThe group $P$ acts on   $G$
by isometries. Therefore, by Lemma 1, one can introduce a Riemannian
structure of the space of orbits
$M_{\bar{p}}$ such that
$$
\bar{\pi}:G\rightarrow M_{\bar{p}}
$$
is an isometry.

\head
3. Curvature of spaces $M_{\bar p}$
\endhead

'In this chapter we will find
conditions on $\bar p$ under that sectional curvature of $M_{\bar p}$
is positive.

'The following Lemma 7 was proved in \cite{E2} but for the sake of
completeness of explanation
we give it's proof.

\proclaim{Lemma 7}
Let $G$ be a compact Lie group with two-sided invariant metric
$\langle\ ,\ \rangle_0$ and  ${\bold t}\subset {\bold g}$
be a maximal commutative subgroup in  the tangent algebra to $G$. Let
$H\in {\bold t}$. Put  $M=Ad(G)A$ where $A\in {\bold g}$.
Let consider a function
$$
f_H:M\rightarrow R:X\mapsto \langle H,X\rangle_0.
$$
'Then extremal values of $f$ are attained on
$M\cap {\bold t}$.
\endproclaim

{\bf "Proof of Lemma 7.}

Firstly, let assume that an element $H$ is regular, i.e., it is contained
only in one maximal commutative subalgebra.

Let  $X\in M$ be a critical point of  $f_H$. ‡That means that
$d_Xf_H=\langle H,- \rangle_0=0$.  Since ' $M=Ad(G)X$, we have
$T_XM=
ad({\bold g})X$.
'Therefore,
$$
\langle ad({\bold g})X,H \rangle_0=0
$$
and
$$
\langle [Z,X],H\rangle_0=\langle Z,[X,H]\rangle_0=0
$$
for every $Z\in {\bold g}$. Hence, ‡
$[X,H]=0$  and, by regularity of $H$, we conclude
that  $X\in {\bold t}$.

"Let assume now that $H$ is a singular element of ${\bold t}$. Let  $X$
an extremal point and  $X\in {\bold g}\setminus {\bold t}$.
By small perturbation of $H$ we obtain a situation when
$H$ be a regular element and $X$ still lies  in
${\bold g}\setminus {\bold t}$.

Thus we arrive at a contradiction with statement proved before.

Lemma 7 is proved.

\proclaim{Lemma 8}
Let $F$ be a subalgebra of ${\bold su}(4)$ with dimension  10
and ${\bold t}$ be a maximal commutative subalgebra formed by diagonal
matrices. Let assume that  $H\in {\bold t}$ and $\langle H,F\rangle_0=0$.
Then up to permutation there exist only two possibilities :
$$
H=i\cdot t\cdot diag(1,1,-1,-1)
$$
or
$$
H=i\cdot t\cdot diag(1,1,1,-3)
$$
where $t\in R$.
\endproclaim

{\bf "Proof of Lemma 8.}

 Let consider a transformation
$$
ad(H):{\bold su}(4)\rightarrow {\bold su}(4):X\mapsto [H,X].
$$
'Take  $X,Y\in F$. 'Then  $[X,Y]\in F$, i.e.,
$\langle H,[X,Y] \rangle_0=0$. 'Therefore,
$$
\langle [H,Y],X\rangle_0=\langle H,[X,Y]\rangle_0=0.
$$
ˆThus,
$$
\langle F,ad\ H(F)\rangle_0=0
$$
ŠMoreover,
$\langle [H,X],H\rangle_0=\langle
[H,H],X\rangle_0=0$ for every $X\in F$,
i.e.,
$$
\langle ad\ H(F),H\rangle_0=0.
$$
'Then
$$
dim(ad\ H(F))\leq dim\ \ {\bold su}(4)-dimF-1=4.
$$
‡That means that
$$
dim(Ker(ad\ H)\cap F)\geq 10-4=6.
$$
' Since $H\in Ker(ad\ H)$ and  $H$ does not lie in $F$, we have
$$
dim(Ker(ad\ H))\geq dim(Ker(ad\ H)\cap F)+1\geq 7.
$$

 Let consider  the $ad\ H$-invariant root decomposition
$$
{\bold su}(4)=V_0\oplus \bigoplus_{\overset{i,j=1}\to{i<j}}^{4}V_{i,j}
$$
where $V_0={\bold t},dimV_{i,j}=2$. 'Then
$$
ad\ H(V_0)=0,
$$
$$
ad\ H(V_{i,j})=\theta_{i,j}(H)V_{i,j}
$$
where $\theta_{i,j}$ are roots of ${\bold su}(4)$, i.e.,
$\theta_{i,j}(i\cdot diag(x_1,x_2,x_3,x_4))=x_i-x_j$.
By estimating of the dimension of the kernel of $adH$, we derive that at least
two roots vanish at $H$.

Lemma 8 is proved.

\proclaim{Lemma 9}
Let for some $m\in M_{\bar{p}}$ there exists a two-dimensional plane
$\sigma \subset T_mM_{\bar{p}}$ such that $K(\sigma)=0$. 'Then there exist
$g\in G,\ X,Y\in {\bold  u}(5)$, such that   $X,Y$ are linearly independent,
$K(X,Y)=0$ and $X,Y$ are orthogonal to a subspace
$$
D_g=\{Ad(g^{-1})\cdot i\cdot t\cdot diag(p_1,p_2,p_3,p_4,p_5)-i\cdot
s\cdot diag(1,1,1,1,0)-\left( \matrix A & 0 \\  0 & 0
\endmatrix \right) |
$$
$$
| t,s\in R, A\in {\bold sp}(2) \}.
$$
\endproclaim

{\bf "Proof of Lemma 9.}

 Let consider the  Riemannian submersion
$$
\bar{\pi}:G\rightarrow M_{\bar{p}}.
$$
"Put $\bar{\pi}(g)=m$  where $g\in G$.
'Then $T_gG=V\oplus H$ where  $V$ is a vertical subspace , $H$
is a horizontal subspace, and $d_g\pi |_H$ is an isometry.
We have
$$
V=T_gP=\{d_eR_g(X)-d_eL_g(Y)|(X,Y)\in T_eP \},
$$
where  $R_g$ and $L_g$ are right and left shift by  $g$ , respectively.
'Then  $H=V^{\perp}$ is an orthogonal complement.

‡Consequently, there exists $\sigma^{*}\in H$ such that
$d_g\bar{\pi}(\sigma^{*})=\sigma$. By Lemmas 2 and 3 ,
$0\leq K(\sigma^{*}) \leq K(\sigma )$ and  therefore  $K(\sigma^*)=0$.
The left shift   $L_{g^{-1}}=(L_g)^{-1}$ is an isometry.
Let $d_gL_{g^{-1}}(\sigma^{*})=Span(X,Y)$ where $X,Y\in T_eG={\bold u}(5)$.
'Then $K(X,Y)=0$ and $\langle V,\sigma^{*} \rangle=\langle
d_gL_{g^{-1}}(V),Span(X,Y) \rangle = 0$.

' Thus, $X,Y$ are orthogonal to a subspace
$$
D_g=d_gL_{g^{-1}}(V)=
$$
$$
=\{ (d_gL_g)^{-1}(d_gR_g)(X)-Y|(X,Y)\in T_eP \}=
$$
$$
=\{ d_e(L_{g^{-1}}\circ R_g)(X)-Y|(X,Y)\in T_eP\}=
$$
$$
=\{ Ad( g^{-1} )X-Y | X\in T_e S^1,Y\in T_e(Sp(2)\times S^1)\}=
$$
$$
=\{ Ad(g^{-1})\cdot
i\cdot t\cdot diag(p_1,p_2,p_3,p_4,p_5)-i\cdot
s\cdot diag(1,1,1,1,0)-\left( \matrix A & 0 \\  0 & 0
\endmatrix \right) |
$$
$$
| t,s\in R, A\in {\bold sp}(2)\}.
$$
Lemma 9 is proved.

\proclaim{Theorem 1}
Let $\bar{p}=(p_1,p_2,p_3,p_4,p_5)$ satisfies to the following conditions

1) $p_{\sigma (1)}+p_{\sigma (2)}-p_{\sigma (3)}-p_{\sigma (4)}$
is relatively prime with $p_{\sigma (5)}$ ;

2) $p_1,p_2,p_3,p_4,p_5>0$;

3) $p_{\sigma (1)}+p_{\sigma (2)}+p_{\sigma (3)}>p_{\sigma
(4)}+p_{\sigma (5)}$;

4) $p_{\sigma (1)}+p_{\sigma (2)}+p_{\sigma (3)}+p_{\sigma
(4)}>3p_{\sigma (5)}$;

5) $3(p_{\sigma (1)}+p_{\sigma (2)})>p_{\sigma (3)}+p_{\sigma
(4)}+p_{\sigma (5)}$

for every permutation $\sigma \in S_5$.

'Then $M_{\bar{p}}$ is positively curved.
\endproclaim

{\bf "Proof of Theorem 1.}

"
Let assume that the statement of Theorem is not valid. Then, by Lemma 9,
there exist $g\in G$ and $X,Y\in {\bold u}(5)$ such that $X,Y$ are
linearly independent and orthogonal to $D_g$ and $K(X,Y)=0$.
By Lemma 5, we may suppose that  '
$Y\in {\bold k}={\bold u}(4)\oplus {\bold u}(1)$ and
$$
[X_p,Y]=[X_k,Y]=0.
$$
 Let consider two possible cases.

{\bf 'Case 1}: $X\in {\bold k}$. 'Then
$$
X=\left( \matrix X_1 & 0 \\  0 & it \endmatrix
\right),\ Y=\left( \matrix Y_1 & 0 \\  0 & is
\endmatrix \right)
$$
where $X_1,Y_1\in {\bold u}(4),t,s\in R$. "The condition that
$[X,Y]=0$ means that
$[X_1,Y_1]=0$. It follows from orthogonality  to ˆ$D_g$ that
$$
\langle X,\left( \matrix A & 0 \\  0 & 0 \endmatrix
\right) \rangle=\langle Y,\left( \matrix A & 0 \\  0 &
0 \endmatrix \right) \rangle=0, \forall A \in {\bold sp}(2),
$$
$$
\langle X,i\cdot diag(1,1,1,1,0) \rangle=\langle Y,i\cdot
diag(1,1,1,1,0) \rangle=0,
$$
$$
\langle X,Ad(g^{-1})\cdot i\cdot
diag(p_1,p_2,p_3,p_4,p_5)\rangle=\langle Y,Ad(g^{-1})\cdot i\cdot
diag(p_1,p_2,p_3,p_4,p_5)\rangle=0.
$$
Since $X,Y\in {\bold k}$, by Lemma 4
the last equalities also valid for the metric
 $\langle \,\ \rangle_0$ .  Therefore $X_1,Y_1$ are orthogonal to
$i\cdot diag(1,1,1,1)$  that means that $X_1,Y_1 \in {\bold su}(4)$ and
$$
\langle X_1,{\bold sp}(2)
\rangle_0=\langle Y_1,{\bold sp}(2) \rangle_0=0,
$$
where ${\bold sp}(2)$
is standard embedded into ${\bold su}(4)$.

ˆIt is known that  $SU(4)/Sp(2)$
is a symmetric rank one space, diffeomorphic to
 $S^{5}$, with positively curved metric.
Therefore $X_1$ and $Y_1$ are linearly dependent. Consequently,
we may suppose that $X_1=0$ with the same  $Span(X,Y)$.

Hence, we may suppose that
$$
X=i\cdot \left( \matrix 0 & 0 \\  0 & 1 \endmatrix
\right)
$$
and
$\langle X,Ad(g^{-1})\cdot i\cdot
diag(p_1,p_2,p_3,p_4,p_5)\rangle_0=0.
$
 Let consider a function
$$
f_X:Ad(G)\cdot i\cdot diag(p_1,p_2,p_3,p_4,p_5)\rightarrow R:Z\mapsto
\langle Z,X\rangle_0.
$$
'By Lemma 7,  $f_X$ attains it's extremal values at points of
$$
Ad(G)\cdot i\cdot diag(p_1,p_2,p_3,p_4,p_5)\cap \{ diagonal\ \ \
matrices \}.
$$
Since two adjoint diagonal matrices coincide up to permutation of
their elements,
the extremal values of  $f_X$
are contained in  $\{ p_1,p_2,p_3,p_4,p_5 \} \subset (0,\infty)$
and thus we arrive at a contradiction.

{\bf Case' 2}: $X$ does not lie in  ${\bold k}$.
Then '
$$
X_p=\left( \matrix 0 & x \\  -x^* & 0 \endmatrix
\right),\ Y=\left( \matrix Y_1 & 0 \\  0 & it
\endmatrix \right),
$$
where $t\in R,Y_1\in {\bold u}(4),x\in C^4\setminus 0$,
and $[X_p,Y]=0$, i.e.,
$$
[X_p,Y]=\left( \matrix 0 & itx \\  -x^*Y_1
& 0 \endmatrix \right)-\left( \matrix 0 & Y_1x \\
-itx^* & 0 \endmatrix \right)=
$$
$$
=\left( \matrix 0 & itx-Y_1x \\  i
tx^*-x^*Y_1 & 0 \endmatrix \right)=0.
$$
‡
Consequently,
$$
Y_1x=itx.
$$
' Since $\langle Y,i\cdot diag(1,1,1,1,0) \rangle=0$, we have
 $\langle Y,i\cdot diag(1,1,1,1,0) \rangle_0=0$ and, therefore,
$$
Y_1\in {\bold su}(4).
$$
' Since $Y_1x=itx$, there exists  $h_1\in SU(4)$ such that
$h_1Y_1h_1^{-1}=i\cdot diag(s_1,s_2,s_3,t)$ where  $s_1+s_2+s_3+t=0$.
ŽDenote
$$
h=\left( \matrix h_1 & 0 \\  0 & 1 \endmatrix
\right)\in SU(5),
$$
then
$$
Y=h^{-1}\cdot i\cdot diag(s_1,s_2,s_3,t,t) \cdot h.
$$
Let
$$
H_1=i\cdot diag(s_1,s_2,s_3,t)\in {\bold su}(4), H=i\cdot
diag(s_1,s_2,s_3,t,t)\in {\bold u}(5).
$$
"The following condition
$$
\langle Y,\left( \matrix A & 0 \\  0 & 0 \endmatrix
\right) \rangle=\frac{1}{2}\langle Y,\left( \matrix A & 0 \\
 0 & 0 \endmatrix \right) \rangle_0=0
$$
means that
$$
\langle Y_1,A\rangle_0=0
$$
for every $ A\in {\bold sp}(2)$.
By using of two-sided invariance of the metric
ˆ $\langle\ \,\ \rangle_0$,
we obtain :
$$
\langle H_1,h_1{\bold sp}(2)h_1^{-1}\rangle_0=0.
$$
ˆIt immediately follows from Lemma 8 that there are two possible values of
$H_1$ up to permutation of coordinates:
$$
H_1=i\cdot t\cdot diag(1,1,-1,-1)
$$
or
$$
H_1=i\cdot t\cdot diag(1,1,1,-3)
$$
where $ t\in R$.
Therefore we may suppose that $H$ satisfies  one of three possibilities :
$$
H=i\cdot diag(1,1,1,-1,-1),
$$
$$
H=i\cdot diag(1,1,1,1,-3),
$$
$$
H=i\cdot diag(3,3,-1,-1,-1).
$$
ŽIt remains to consider a condition
$$
0=2\langle Y,Ad(g^{-1})\cdot i\cdot diag(p_1,p_2,p_3,p_4,p_5) \rangle=
$$
$$
=\langle Y,Ad(g^{-1})\cdot i\cdot diag(p_1,p_2,p_3,p_4,p_5)
\rangle_0=
$$
$$
=\langle h^{-1}Hh,Ad(g^{-1})\cdot i\cdot
diag(p_1,p_2,p_3,p_4,p_5) \rangle_0=
$$
$$
=\langle H,Ad(g')\cdot i\cdot diag(p_1,p_2,p_3,p_4,p_5) \rangle_0,
$$
where $g'=hg^{-1}\in G$.

 Let consider a function
$$
f_H:Ad(G)\cdot i\cdot diag(p_1,p_2,p_3,p_4,p_5)\rightarrow R:X\mapsto
\langle X,H\rangle_0.
$$
By Lemma 7, it's extremal values are attained at the set
$$
\{
p_{i_1}+p_{i_2}+p_{i_3}-p_{i_4}-p_{i_5},
3(p_{i_1}+p_{i_2})-p_{i_3}-p_{i_4}-p_{i_5},
p_{i_1}+p_{i_2}+p_{i_3}+p_{i_4}-3p_{i_5} |
$$
$$
| \{ i_1,i_2,i_3,i_4,i_5 \}=\{ 1,2,3,4,5 \} \}
$$
which lies in  $(0,\infty)$ by the condition of Theorem.
Thus we arrive at a contradiction.

Theorem 1 is proved.

‹It is easy to see that all conditions of Theorem 1 hold for
$p_1=1,p_2=p_3=p_4=p_5=q^n$ where $q$ is a prime number and
$ n$ is a nonnegative integer.

\head
4. Topology of spaces $M_{\bar p}$
\endhead

ŽLet denote by $\sigma_i(\bar{p})$ the $i$-th elementary symmetric
function of $p_1,p_2,p_3,p_4,p_5$.

\proclaim{Lemma 10}
$(Sp(2)\times S^1)/\pm(E,1)$
is diffeomorphic to $Sp(2)\times S^1$.
\endproclaim

{\bf "Proof of Lemma 10.}

 Let consider the mapping
$$
\phi :Sp(2)\times S^1\rightarrow Sp(2)\times S^1:(A,z)\mapsto (A\cdot
diag(z,z,\bar{z},\bar{z}),z^2).
$$
Put $\phi '(A,z)=(B,w)$. 'Then $z^2=w,z=\pm \sqrt{w}$ that means
$$
(A,z)=\pm (B\cdot
diag(\sqrt{\bar{w}},\sqrt{\bar{w}},\sqrt{w},\sqrt{w}),\sqrt{w}).
$$

' Thus, the mapping  $\phi '$  induces a bijection
$$
\phi :\frac{Sp(2)\times S^1}{\pm (E,1)} \rightarrow Sp(2)\times
S^1
$$
which evidently occurs to be a diffeomorphism.

Lemma 10 is proved.

\proclaim{ 'Theorem 2}
Let
$(p_{\sigma (1)}+p_{\sigma (2)}-p_{\sigma (3)}-p_{\sigma (4)})$ and
$p_{\sigma (5)}$ are relatively prime for every permutation
$\sigma \in S_5.$
'Then the space
$M_{\bar{p}}$ is simply connected.
\endproclaim

{\bf "Proof of Theorem 2.}

 Let consider the fragment of the exact homotopy sequence of the fiber bundle
$\bar{\pi}:G\rightarrow M_{\bar{p}}$ with the fiber
$P=S^1\times (Sp(2)\times
S^1)/\pm (E,1)$:
$$
\pi_1(S^1\times \frac{Sp(2)\times
S^1}{\pm (E,1)})\overset{i_*}\to{\rightarrow}
\pi_1(U(5))\overset{\bar{\pi}_*}\to{\rightarrow}
\pi_1(M_{\bar{p}})\rightarrow 0
$$
where $i$ is an embedding of  $P$ as the fiber over the unit element
$E\in U(5)$. ˆ
Since $\phi$ is a diffeomorphism, we obtain
$$
\pi_1(S^1\times Sp(2)\times
S^1)\overset{j_*}\to{\rightarrow}\pi_1(U(5))\overset{\bar{\pi}_*}\to
{\rightarrow}\pi_1(M_{\bar{p}})\rightarrow 0,
$$
where
$j=i\circ (id\times \phi^{-1})$.
Thus we have the following exact sequence
$$
{\bold Z}\oplus {\bold Z}\overset{j_*}\to{\rightarrow}
{\bold Z}\overset{\bar{\pi}_*}\to{\rightarrow} \pi_1(M_{\bar{p}})\rightarrow 0.
$$
'Let compute $j_*$.
'We take in the group ${\bold Z}\oplus {\bold Z}=\pi_1(S^1\times Sp(2)\times
S^1)$
it's generators $(1,0)$ and  $(0,1)$  which are realized by loops
$$
\xi_1(t)=(e^{2\pi it},(E,1)),\ \
\xi_2(t)=(1,(E,e^{2\pi it})),\ \
0\leq t\leq 1.
$$
'We take by a generator in ${\bold Z}=\pi_1(U(5))$ a homotopy class given by
the
following winding of torus :
$$
\xi (t)=diag(e^{2\pi ix_1t},e^{2\pi ix_2t},e^{2\pi ix_3t},e^{2\pi
ix_4t},e^{2\pi ix_5t}),0\leq t\leq 1,
x_k \in Z, \sum_{k=1}^{5}x_k=1.
$$ 
The loop $\xi_1$ is mapped into
$$
j(\xi_1(t))=i(e^{2\pi it},\pm (E,1))=diag(e^{2\pi ip_1t},e^{2\pi
ip_2t},e^{2\pi ip_3t},e^{2\pi ip_4t},e^{2\pi ip_5t}),\ 0\leq t\leq 1
$$
 and the loop $\xi_2$ is mapped into
$$
j(\xi_2(t))=i(1,\pm (diag(e^{-\pi it},e^{-\pi it},e^{\pi it},e^{\pi
it}),e^{\pi it}))=diag(1,1,e^{2\pi it},e^{2\pi it},1),\
0\leq t\leq 1.
$$
We obtain that
$$
j_*(1,0)=\sigma_1(\bar{p})\cdot 1,\ \
j_*(0,1)=2\cdot 1.
$$
ˆIt follows, in particular,  from conditions of Theorem that
2 and  $\sigma_1(\bar{p})$ are relatively prime and, since that, we conclude
that  $j_*$ is an epimorphism.
It immediately follows from the exact sequence
quoted above that the manifold
$M_{\bar{p}}$ is simply connected.

'Theorem 2 is proved.

ˆPut  $G=U(5)$ and $P=S^1\times (Sp(2)\times S^1)/{\bold Z}_2\subset
G\times G$. We denote by  $M=M_{\bar{p}}=G/P$  the space of orbits, and
denote by $\bar{\pi}:G\rightarrow M$ the principal bundle with the
structure group $P$.
Let $\pi_G:E_G\rightarrow B_G$ and  $\pi_P:E_P\rightarrow B_P$
be universal coverings for groups $G$ and  $P$, respectively, with contractible
covering spaces $E_G$ and $E_P$.
Let consider the following commutative diagram :
$$
\matrix
 & G &
\overset{p_2}\to{\longleftarrow} & E_P\times G &
\overset{p_1}\to{\longrightarrow} & E_P & \\ \pi_G & \downarrow & &
\downarrow & & \downarrow & \pi_P \\ & M &
\overset{\bar{p}_2}\to{\longleftarrow} & G//P &
\overset{\bar{p}_1}\to{\longrightarrow} & B_P &
\endmatrix
$$
‡Here $p_1$ and  $p_2$ are natural projections onto first and second
factors, $G//P$ is the space of orbits of the natural action of $P$ on
$E_P\times G$. Since the  fiber of $\bar{p}_2$ is diffeomorphic to contractible
space $E_P$, $\bar{p}_2^*$ maps $H^*(M)$ isomorphically onto
$H^*(G//P)$.   Let consider the spectral sequence of the fiber bundle
$p=\bar{p}_1:G//P\rightarrow B_P$ with the fiber $G$. ' Since $H^*(G)$ is
torsion free, the initial term $E_2=H^*(B_P)\otimes H^*(G)$.
The term  —$E_\infty$ is attached to $H^*(M)$. Let compute differentials of
this
spectral sequence.

We consider the diagram:
$$
\matrix
 G//P= & (E_{G^2}\times G)/P &
\overset{\hat{\rho}}\to{\longrightarrow} & (E_{G^2}\times G)/G^2 &
\overset{f}\to{\longleftarrow} & B_G & =E_{G^2}/\delta G \\
& p\downarrow & & \downarrow p' & & \downarrow \triangle & \\ &
 B_P & \overset{\rho}\to{\longrightarrow} & B_{G^2} &
\overset{id}\to{\longleftarrow} & B_{G^2} &
\endmatrix
$$
‡Here we put $E_P=E_{G^2}$, $B_P=E_{G^2}/P$, and
$B_{G^2}=E_{G^2}/G^2$. We denote here by $\rho :B_P\rightarrow B_{G^2}$
and $\triangle:B_G\rightarrow B_{G^2}$ natural projections. We also denote by
$\delta:G\rightarrow G^2$ the diagonal embedding and denote by
$\hat{\rho}$ the fibered mapping whose
restrictions onto fibers are homeomorphisms.
The mapping $f:(\delta G)e\mapsto G^2(e,1)$ is an isomorphism of
fiber bundles.

'Let compute differentials of the spectral sequence of the fiber bundle
$\triangle$. ŠWe identify the ring of cohomologies $H^*(G)$
with an interior algebra with generators $z_1,z_3,\ldots,z_9$.
Then we put
$H^*(B_G)={\bold Z}[\bar{z}_1,\bar{z}_3,\ldots,\bar{z}_9]$ where $\bar{z}_i$
s the image of $z_i$ under transgression of the fiber bundle $\pi_G$.
Since $B_{G^2}=B_G\times B_G$,
$$
H^*(B_{G^2})=H^*(B_G)\otimes
H^*(B_G)=
{\bold Z}[\bar{x}_1,\bar{y}_1,\bar{x}_3,\bar{y}_3,\ldots,\bar{x}_9,\bar{y}_9]
$$
where  $\bar{x}_i=\bar{z}_i\otimes 1, \bar{y}_i=1\otimes \bar{z}_i$.

 The initial term $E_2$ is isomorphic to $H^*(B_{G^2})\otimes H^*(G)$.
Let denote by  Ž
$k_i:H^*(B_{G^2})\rightarrow E_i^{*,0}$ the natural projection. Then it is
well-known that
$$
\triangle^*=k_{\infty}:H^*(B_{G^2})\rightarrow
E_{\infty}^{*,0}\subset H^*B_G.
$$

\proclaim{Lemma 11}
$$
1)\ d_j(1\otimes z_i)=0,\ \ \ \  j\leq i,\ i=3,5,7 ;
$$
$$
2)\ d_{i+1}(1\otimes z_i)=\pm k_{i+1}(\bar{x}_i-\bar{y}_i),\
i=1,3,5,7.
$$
\endproclaim

{\bf "Proof of Lemma 11.}

  Let consider the term  $E_2$ of the spectral sequence of the fiber bundle
$B_G\rightarrow B_{G^2}.$
$$
\matrix z_7 & 0 & \ast & 0 &
\ast & 0 & \ast  \\  z_1z_5 & 0 &
z_1z_5\bar{x}_1, z_1z_5\bar{y}_1 & 0 & \ast & 0 & \ast
\\  z_5 & 0 & \ast & 0 & \ast & 0 & \ast \\
z_1z_3 & 0 & z_1z_3\bar{x}_1,z_1z_3\bar{y}_1 & 0 & z_1z_3\otimes
E_2^{4,0} & 0 & \ast \\  z_3 & 0 &
z_3\bar{x}_1,z_3\bar{y}_1 & 0 & \ast & 0 & \ast \\  0
& 0 & 0 & 0 & 0 & 0 & 0 \\  z_1 & 0 &
z_1\bar{x}_1,z_1\bar{y}_1 & 0 &
z_1\bar{x}_1^2,z_1\bar{x}_1\bar{y}_1,z_1\bar{y}_1^2 & 0 & \ast
\\ & & & & z_1\bar{x}_3,z_1\bar{y}_3 & &  \\  1 & 0 &
\bar{x}_1,\bar{y}_1 & 0 & \bar{x}_1^2,\bar{x}_1\bar{y}_1,\bar{y}_1^2 &
0 & \bar{x}_1^3,\bar{x}_1^2\bar{y}_1,\bar{x}_1\bar{y}_1^2,\bar{y}_1^3
\\ & & & & \bar{x}_3,\bar{y}_3 & &
\bar{x}_1\bar{x}_3,\bar{x}_1\bar{y}_3,\bar{y}_1\bar{x}_3,\bar{y}_1\bar{y}_3
\\ & & & & & & \bar{x}_5,\bar{y}_5 \\  \endmatrix
$$

ˆWe have
$$
\triangle^*(1\otimes u)=\triangle^*(u\otimes 1)=u.
$$
Ÿfor every $u \in H^*(B_G)$.
Since the kernel of  $\triangle^2$ coincides with $d_2({\bold Z}(z_1))$,
$d_2({\bold Z}(z_1))={\bold Z}(\bar{x}_1-\bar{y}_1)$. ‡We conclude that
$$
d_2(z_1)=\pm (\bar{x}_1-\bar{y}_1)=\pm k_2(\bar{x}_1-\bar{y}_1).
$$
' Thus, (2) is proved for $i+1$.

" Then we have
$$
d_2(z_1z_3\bar{x}_1)=z_3\bar{x}_1^2-z_3\bar{x}_1\bar{y}_1,
$$
$$
d_2(z_1z_3\bar{y}_1)=z_3\bar{x}_1\bar{y}_1-z_3\bar{y}_1^2.
$$
'Therefore, $Ker(d_2^{2,4})=0$ and we get that $d_2^{0,5}=0$. €
Analogously, we conclude that $Ker(d_2^{2,6})=0, Ker(d_2^{4,4})=0$ and,
therefore,
$d_2^{0,7}=d_4^{0,7}=0$.
'Triviality of other differentials follows from
dimensional reasons.
Thus, it is left to prove (2) for  $i=3,5,7$.
‹One can easily see that
$$
Ker\triangle^4={\bold Z}(\bar{x}_3-\bar{y}_3)\oplus
{\bold Z}(\bar{x}_1^2-\bar{x}_1\bar{y}_1,\bar{x}_1\bar{y}_1-\bar{y}_1^2),
$$
and, from other side,
$$
Ker\triangle^4=Im(d_2^{2,1})\oplus
Im(d_4^{0,3}).
$$
"Taking into account that
$$
Im(d_2^{2,1})={\bold Z}(\bar{x}_1^2-\bar{x}_1\bar{y}_1,
\bar{x}_1\bar{y}_1-\bar{y}_1^2),
$$
we derive that
$$
d_4(z_3)=\pm k_4(\bar{x}_3-\bar{y}_3).
$$
'In the same manner we obtain that
$$
Ker\triangle^6={\bold Z}(\bar{x}_5-\bar{y}_5)\oplus
{\bold Z}(\bar{x}_1^3-\bar{x}_1^2\bar{y}_1,\bar{x}_1^2\bar{y}_1-
\bar{x}_1\bar{y}_1^2,\bar{x}_1\bar{y}_1^2-\bar{y}_1^3)\oplus
$$
$$
\oplus
{\bold Z}(\bar{x}_1\bar{x}_3-\bar{y}_1\bar{x}_3,
\bar{x}_1\bar{y}_3-\bar{y}_1\bar{y}_3,
\bar{x}_1\bar{x}_3-\bar{x}_1\bar{y}_3),
$$
and
$$
Ker\triangle^6=Im(d_2^{4,1})\oplus Im(d_4^{2,4})\oplus Im(d_6^{0,5}).
$$
"Taking into account that we proved before that  first two summands
from  the last expression coincides respectively with the same from the
preceding one we get
$$
d_6(z_5)=\pm k_6(\bar{x}_5-\bar{y}_5).
$$
' Analogously one can prove that $d_8(z_7)=\pm k_8(\bar{x}_7-\bar{y}_7).$

Lemma 11 is proved.

\proclaim{Lemma 12}
Let $d_j, j \geq 1$ be differentials in the spectral sequence
of the fiber bundle  $p:G//P\rightarrow B_P$. Then

1)\ $d_j(1\otimes z_i)=0,\ \ \ j\leq i,\ i=3,5,7$ ;

2)\ $d_{i+1}(1\otimes z_i)=\pm k_{i+1}\rho^* (\bar{x}_i-\bar{y}_i),\
i=1,3,5,7$

\noindent
where $\rho :B_P\rightarrow B_{G^2}$ is induced by the embedding
$P\subset G^2$.
\endproclaim

{\bf "Proof of Lemma 12.}

 Let consider the second diagram. The fibered mapping
$(\hat{\rho},\rho)$ generates the homomorphism $\rho^{\sharp}$, of spectral
sequences, and, moreover,  $\rho_2^{\sharp}=\rho^*\otimes
i:H^*B_{G^2}\otimes H^*G\rightarrow H^*B_P\otimes H^*G$ where
$i$ is an isomorphism.
We put $i(1\otimes z_i)=1\otimes z_i.$ ' Then the following identities
$$
d_j(1\otimes z_i)=\rho^{\sharp}(d'_j(1\otimes z_i))=\rho^
{\sharp}(0)=0,j\leq i,
$$
$$
d_{i+1}(1\otimes
z_i)=\rho^{\sharp}(d'_{i+1}(1\otimes z_i))= \pm \rho^{\sharp}
(k'_{i+1}(\bar{x}_i-\bar{y}_i))=\pm
k_{i+1}(\rho^*(\bar{x}_i-\bar{y}_i)
$$
hold.

Lemma 12 is proved.


Let $G$ be a Lie group and $T^n$ be a maximal torus in
$G$ where
$i:T^n\rightarrow G$ is an embedding and $j:B_{T^n}\rightarrow B_G$
is a natural projection. We denote by
$a_1,\ldots ,a_n$ generators of
$H^1T^n$. Then we have  $H^*B_{T^n}={\bold Z}[\bar{a}_1,\ldots ,\bar{a}_n]$.
Let denote by
$I_G$ the algebra of polynoms in  $H^*B_{T^n}$ that are invariant under the
action of the Weyl group $W(G)$.

\proclaim{Borel Theorem}
 (\cite{Bo})
Let …$H^*G$ and  $H^*(G/T^n)$ are torsion free. Then
$j^*:H^*B_G\rightarrow H^*B_{T^n}$ is a monomorphism and it's image coincides
with $I_G$.
\endproclaim

Š Borel proved (\cite{Bo}) that conditions of this theorem hold for
every classic group. We have $G=U(5)$, $\bar{z}_1=\sigma_1(\bar{d}_1,\ldots
,\bar{d}_5),z_3=\sigma_2(\bar{d}_1,\ldots ,\bar{d}_5),\ldots
,z_9=\sigma_5(\bar{d}_1, \ldots ,\bar{d}_5)$ where   $d_1,\ldots ,d_5$
are cocycles which are adjoint to cycles  $D_1,\ldots ,D_5$ that are defined by
$D_i(t)=diag(1,\ldots,e^{2\pi it},\ldots ,1),0\leq t\leq 1$.

 Let consider $P\subset G^2$:
$$
P=S^1\times \frac{Sp(2)\times S^1}{{\bold Z}_2}\simeq S^1\times S^1\times
Sp(2),
$$
ŠThe ring of cohomologies $H^*P$ is torsion free.
Next, let  $T^3$ be a maximal torus in $Sp(2)\times S^1$. 'Then
$T^3/{\bold Z}_2$ is a maximal torus in $(Sp(2)\times S^1)/{\bold Z}_2$, i.e.,
$$
\frac{S^1\times
Sp(2)}{{\bold Z}_2}/\frac{T^3}{{\bold Z}_2}\simeq \frac{S^1\times Sp(2)}{T^3}.
$$
One can see now that for $P$ the conditions of Borel Theorem hold.

Let $T$ be a torus in $G^2$ and $S$ be a torus in $P$. Denote by
$i:S\rightarrow T$ an embedding and  denote by
$j:B_S\rightarrow B_T$ a natural projection.
Let consider the following diagram
$$
\matrix
 & H^*B_S &
\overset{j^*}\to{\longleftarrow} & H^*B_T & \\ j_S^* & \uparrow & &
\uparrow & j_T^* \\ & H^*B_P & \overset{\rho^*}\to{\longleftarrow} &
H^*B_{G^2} &
\endmatrix
$$
Let choose a basis  $A_1,\ldots ,A_5,B_1,\ldots ,B_5$
of cycles in  $H_1(T)$ as follows :
$$
A_i(t)=(1,diag(1,\ldots,e^{2\pi
it},\ldots ,1)),0\leq t\leq 1,
$$
$$
B_i(t)=(diag(1,\ldots,e^{2\pi
it},\ldots ,1),1),0\leq t\leq 1,
$$
and denote by
$a_1,\ldots ,a_5,b_1,\ldots ,b_5$
cocycles ($\in H^1(T)$) that are adjoint to elements of this basis.
Let choose a basis $C_1,\ldots ,C_4$
of cycles in  $H_1(S)$ as follows :
$$
C_1(t)=(e^{2\pi it},\pm (E,1)),
$$
$$
C_2(t)=(1,\pm (diag(e^{\pi it},e^{\pi it},e^{-\pi it},e^{-\pi
it}),e^{\pi it})),
$$
$$ C_3(t)=(1,\pm (diag(e^{2\pi it},1,e^{-2\pi
it},1),1)),
$$
$$
C_4(t)=(1,\pm (diag(1,e^{2\pi it},1,e^{-2\pi
it}),1)),
$$
$$
0\leq t\leq 1,
$$
and denote by $c_1,c_2,c_3,c_4$ cocycles ($H^1(S)$)
that are adjoint to elements
of this basis. Then we have
$$
i_*(C_1)=p_1B_1+p_2B_2+p_3B_3+p_4B_4+p_5B_5,
$$
$$
i_*(C_2)=A_1+A_2, \ \
i_*(C_3)=A_1-A_3, \ \
i_*(C_4)=A_2-A_4.
$$
'Consequently,
$$
i^*(a_1)=c_2+c_3, \ \
i^*(a_2)=c_2+c_4, \ \
i^*(a_3)=-c_3, \ \
i^*(a_4)=-c_4, \ \
i^*(a_5)=0,
$$
$$
i^*(b_1)=p_1c_1, \ \
i^*(b_2)=p_2c_1, \ \
i^*(b_3)=p_3c_1, \ \
i^*(b_4)=p_4c_1, \ \
i^*(b_5)=p_5c_1.
$$
'Since transgression is natural, we have
$$
j^*(\bar{a}_i)=\overline{i^*(a_i)},\ \
j^*(\bar{b}_i)=\overline{i^*(b_i)}.
$$
By the diagram given above,
$\rho^*$ is the restriction of $j^*$ onto  $I_{G^2}$.

ŒWe identify $H^*B_{G^2}$
with the subalgebra, of $H^*B_T$,
generated by
$$
\sigma_i(\bar{a}_1,\ldots ,\bar{a}_5), \ \ \
\sigma_i(\bar{b}_1,\ldots ,\bar{b}_5),  \ \ \
i=1,2,\ldots,5.
$$
ŠWe identify the ring of cohomologies $H^*B_P$ with subalgebra,
 in $H^*B_S$, that is invariant under $W(P)$. ‡ Notice that in our
definitions
$$
a_i=1\otimes
d_i,\ b_i=d_i\otimes 1.
$$
'Let compute  $W(P)$. Elements of  $W(P)$
are induced by elements of  $W(S^1\times S^1\times Sp(2))$.
'Hence, generators $\phi_1,\phi_2,\phi_3$, of  $W(P)$,
act on homologies of $S$ as follows
$$
\matrix \phi_1: & C_1\mapsto
C_1 & \phi_2: & C_1\mapsto C_1 & \phi_3: & C_1\mapsto C_1 \\ &
C_2\mapsto C_2-C_3 & & C_2\mapsto C_2-C_4 & & C_2\mapsto C_2 \\ &
C_3\mapsto -C_3 & & C_3\mapsto C_3 & & C_3\mapsto C_4 \\ & C_4\mapsto
C_4 & & C_4\mapsto -C_4 & & C_4\mapsto C_3, \endmatrix
$$
that means that their action on cohomologies is given by
$$ \matrix \phi_1: &
c_1\mapsto c_1 & \phi_2: & c_1\mapsto c_1 & \phi_3: & c_1\mapsto c_1 \\
& c_2\mapsto c_2 & & c_2\mapsto c_2 & & c_2\mapsto c_2 \\ & c_3\mapsto
-c_2-c_3 & & c_3\mapsto c_3 & & c_3\mapsto c_4 \\ & c_4\mapsto c_4 & &
c_4\mapsto -c_2-c_4 & & c_4\mapsto c_3. \endmatrix
$$
ˆThus $H^*B_P$ is a subalgebra of ${\bold Z}[\bar{c}_1,\bar{c}_2]$
that is invariant under $W(P)$. Let find multiplicative generators of
$H^*B_P$.

\proclaim{Lemma 13}
Let denote
$$
\bar{f}=(\bar{c}_3^2+\bar{c}_4^2)+\bar{c}_2(\bar{c}_3+\bar{c}_4),
$$
$$
\bar{g}=\bar{c}_3^2\bar{c}_4^2+\bar{c}_2\bar{c}_3\bar{c}_4(\bar{c}_3+
\bar{c}_4)+
\bar{c}_2^2\bar{c}_3\bar{c}_4.
$$

'Then $H^*B_P={\bold Z}[\bar{c}_1,\bar{c}_2,\bar{f},\bar{g}].$
\endproclaim

{\bf "Proof of Lemma  13.}

 Let consider the natural embedding
$H^*B_S={\bold Z}[\bar{c}_1,\bar{c}_2,\bar{c}_3,\bar{c}_4] \subset
{\bold R}[\bar{c}_1,\bar{c}_2,\bar{c}_3,\bar{c}_4]$. '
We denote by $A_{\bold R}$ he subalgebra of   
${\bold R}[\bar{c}_1,\bar{c}_2,\bar{c}_3,\bar{c}_4]$
that is invariant under
$W(P)$. 'Then we have $H^*B_P=A_{\bold Z}=A_{\bold R} \cap
{\bold Z}[\bar{c}_1,\bar{c}_2,\bar{c}_3,\bar{c}_4].$

‡ We define an isomorphism $\tau :{\bold R}[x_1,x_2,x_3,x_4]\rightarrow
{\bold R}[\bar{c}_1,\bar{c}_2,\bar{c}_3,\bar{c}_4]$ as follows:
$$
\matrix
 x_1 \mapsto \bar{c}_1 \\ x_2 \mapsto \bar{c}_2 \\ x_3
\mapsto \bar{c}_3+\frac{1}{2} \bar{c}_2 \\ x_4 \mapsto
\bar{c}_4+\frac{1}{2}\bar{c}_2.
\endmatrix
$$
ƒThe  $W(P)$ is conjugated by  $\tau$ ­ 
to the group $W'$ that acts on  ${\bold R}[x_1,x_2,x_3,x_4]$
and generated by the following elements ;
$$
\matrix
 \phi_1':  & x_1 \mapsto x_1 & \phi_2': & x_1
\mapsto x_1 & \phi_3': & x_1 \mapsto x_1 \\ & x_2 \mapsto x_2 & & x_2
\mapsto x_2 & & x_2 \mapsto x_2 \\ & x_3 \mapsto x_4 & & x_3 \mapsto
-x_3 & & x_3 \mapsto x_3 \\ & x_4 \mapsto x_3 & & x_4 \mapsto x_4 & &
x_4 \mapsto -x_4.
\endmatrix
$$

' Thus, we obtain that  $W'$ is the Weyl group of group
$S^1\times S^1\times Sp(2)$ and it is absolutely evident that
the subalgebra of
${\bold R}[x_1,x_2,x_3,x_4]$ that is invariant under $W'$ coincides with
$A'_{\bold R}={\bold R}[x_1,x_2,x_3^2+x_4^2,x_3^2x_4^2]$.
Hence,
$A_{\bold R}={\bold R}[\tau (x_1),\tau (x_2),\tau (x_3^2+x_4^2),\tau
(x_3^2x_4^2)]$.
Nondifficult computations give
$$
\tau (x_1)=\bar{c}_1,
\tau (x_2)=\bar{c}_2,
$$
$$
\tau
(x_3^2+x_4^2)=(\bar{c}_3^2+\bar{c}_4^2)+\bar{c}_2(\bar{c}_3+\bar{c}_4)+
\frac{1}{2}\bar{c}_2^2,
$$
$$
\tau
(x_3^2x_4^2)=\bar{c}_3^2\bar{c}_4^2+\bar{c}_2\bar{c}_3\bar{c}_4(\bar{c}_3+
\bar{c}_4)+\bar{c}_2^2\bar{c}_3\bar{c}_4+
+\frac{1}{4}\bar{c}_2^2((\bar{c}_3^2+
\bar{c}_4^2)+\bar{c}_2(\bar{c}_3+\bar{c}_4)+\frac{1}{4}\bar{c}_2^2).
$$
' Thus, we have
$
A_{\bold R}={\bold
R}[\bar{c}_1,\bar{c}_2,(\bar{c}_3^2+\bar{c}_4^2)+
\bar{c}_2(\bar{c}_3+\bar{c}_4),
\bar{c}_3^2\bar{c}_4^2+\bar{c}_2\bar{c}_3\bar{c}_4(\bar{c}_3+
\bar{c}_4)+\bar{c}_2^2\bar{c}_3\bar{c}_4]=
{\bold R}[\bar{c}_1,\bar{c}_2,\bar{f},\bar{g}].
$
' Since generators of  $A_{\bold R}$ lie in
${\bold Z}[\bar{c}_1,\bar{c}_2,\bar{c}_3,\bar{c}_4]$, one get
$$
A_{\bold Z}={\bold Z}[\bar{c}_1,\bar{c}_2,\bar{f},\bar{g}].
$$

Lemma 13 is proved.

\proclaim{Theorem 3}
The space $M_{\bar{p}}$ has the following groups of cohomologies
$$
1) H^i=\cases
 {\bold Z}, & for \ \ \ \ $i=0,2,4,9,11,13,$ \cr
 0, & for \ \ \ \ $i=1,3,5,7,10,12$;
\endcases
$$
2) the groups  $H^6(M_{\bar{p}})$ and
$H^8(M_{\bar{p}})$ are finite and their orders are equal to
$$
r=|\sigma_1^3-4\sigma_1\sigma_2+8\sigma_3|.
$$
\endproclaim

{\bf "Proof of Theorem 3.}

ŽWe denote $\sigma_i=\sigma_i(p_1,\ldots,p_5)$.
Let consider the term  $E_2$
of the spectral sequence of the fiber map
$G//P\rightarrow B_P$.
$$
\matrix z_7 & 0 & \ast & 0 &
\ast & 0 & \ast & 0 & \ast \\  z_1z_5 & 0 & z_1z_5\bar{c}_1,
z_1z_5\bar{c}_2 & 0 & \ast & 0 & \ast & 0 & \ast \\  z_5 & 0 &
z_5\bar{c}_1, z_5\bar{c}_2 & 0 & \ast
& 0 & \ast & 0 & \ast \\
 z_1z_3 & 0 &
z_1z_3\bar{c}_1,z_1z_3\bar{c}_2 & 0 & z_1z_3\bar{c}_1^2,
z_1z_3\bar{c}_1\bar{c}_2 & 0 & \ast & 0 & \ast \\ & & & &
z_1z_3\bar{c}_2^2, z_1z_3\bar{f} & & & & \\  z_3 & 0 &
z_3\bar{c}_1,z_3\bar{c}_2 & 0 & z_3\bar{c}_1^2
z_3\bar{c}_1\bar{c}_2 & 0 & \ast & 0 & \ast \\ & & &
&z_3\bar{c}_2^2, z_3\bar{f} & & & & \\  0 & 0 & 0 & 0 & 0 & 0 & 0
& 0 & 0 \\  z_1 & 0 & z_1\bar{c}_1,z_1\bar{c}_2 & 0 &
z_1\bar{c}_1^2,z_1\bar{c}_1\bar{c}_2 & 0 & z_1\bar{c}_1^3,
z_1\bar{c}_1^2\bar{c}_2, z_1\bar{c}_1\bar{c}_2^2 & 0 & \ast \\ & & & &
z_1\bar{c}_2^2,z_1\bar{f} & & z_1\bar{c}_2^3, z_1\bar{c}_1\bar{f},
z_1\bar{c}_2\bar{f} & & \\  1 & 0 & \bar{c}_1,\bar{c}_2 & 0 &
\bar{c}_1^2,\bar{c}_1\bar{c}_2 & 0 &
\bar{c}_1^3,\bar{c}_1^2\bar{c}_2,\bar{c}_1\bar{c}_2^2 & 0 &
\bar{c}_1^3\bar{c}_2, \bar{c}_1^2\bar{c}_2^2,
\bar{c}_1\bar{c}_2^3 \\ & & & & \bar{c}_2^2,\bar{f} & &
\bar{c}_2^3,\bar{c}_1\bar{f},\bar{c}_2\bar{f} & & \bar{c}_2^4,
\bar{c}_1^2\bar{f}, \bar{c}_1\bar{c}_2\bar{f} \\  & & & & & & & &
\bar{c}_2^2\bar{f},\bar{f}^2,\bar{g},\bar{c}_1^4 \\  \endmatrix
$$

ˆWe have
$$
d_2z_1=\pm
\rho^*(\bar{x}_1-\bar{y}_1)=
\rho^*(\sigma_1(\bar{b}_1,\ldots,\bar{b}_5)-
\sigma_1(\bar{a}_1,\ldots,\bar{a}_5))=
$$
$$
=\sigma_1\cdot \bar{c}_1-2\cdot \bar{c}_2.
$$
Take $n\in {\bold Z}$ such that $\sigma_1+2n=1$. 'Then
$$
\frac{{\bold Z}(\bar{c}_1,\bar{c}_2)}{{\bold Z}(\sigma_1\cdot \bar{c}_1-2\cdot
\bar{c}_2)}={\bold Z}(n\cdot \bar{c}_1 +\bar{c}_2).
$$
ŽDenote by $F_2$ the image of $d_2^{2,1}$. Notice that $F_2$ is generated
by elements
$$
d_2(z_1\bar{c}_1)=\sigma_1\bar{c}_1^2-2\bar{c}_1\bar{c}_2,
d_2(z_1\bar{c}_2)=\sigma_1\bar{c}_1\bar{c}_2-2\bar{c}_2^2.
$$
'Then we have
${\bold Z}(\bar{c}_1^2,\bar{c}_1\bar{c}_2,\bar{c}_2^2,\bar{f})/F_2=
{\bold Z}((n-n^2)\bar{c}_1^2+\bar{c}_2^2,\bar{f})$.
ŽDenote by  $F_4$ the image of
$d_2^{4,1}$. We see that  $F_4$ is generated by elements
$$
d_2(z_1\bar{c}_1^2)=\sigma_1\cdot \bar{c}_1^3-2\cdot
\bar{c}_1^2\bar{c}_2 , \ \ \ \  d_2(z_1\bar{c}_1\bar{c}_2)=\sigma_1\cdot
\bar{c}_1^2\bar{c}_2-2\cdot \bar{c}_1\bar{c}_2^2,
$$
$$
d_2(z_1\bar{c}_2^2)=\sigma_1\cdot \bar{c}_1\bar{c}_2^2-2\cdot
\bar{c}_2^3, \ \ \ \ d_2(z_1\bar{f})=\sigma_1\cdot \bar{c}_1\bar{f}-2\cdot
\bar{c}_2\bar{f}
$$
 Finally, denote by $F_6$ the image of
$d_2^{6,1}$ and notice that it is generated by elements
$$
\matrix
d_2(z_1\bar{c}_1^3)=\sigma_1\bar{c}_1^4-2\bar{c}_1^3\bar{c}_2, &
d_2(z_1\bar{c}_2^3)=\sigma_1\bar{c}_1\bar{c}_2^3-2\bar{c}_2^4, \\
d_2(z_1\bar{c}_1^2\bar{c}_2)=
\sigma_1\bar{c}_1^3\bar{c}_2-2\bar{c}_1^2\bar{c}_2^2,
&
d_2(z_1\bar{c}_1\bar{f})=
\sigma_1\bar{c}_1^2\bar{f}-2\bar{c}_1\bar{c}_2\bar{f},
\\
d_2(z_1\bar{c}_1\bar{c}_2^2)=
\sigma_1\bar{c}_1^2\bar{c}_2^2-2\bar{c}_1\bar{c}_2^3, &
d_2(z_1\bar{c}_2\bar{f})=
\sigma_1\bar{c}_1\bar{c}_2\bar{f}-2\bar{c}_2^2\bar{f}.
\endmatrix
$$

Let proceed to the term $E_3=E_4$.
$$
\matrix z_7 & 0 & \ast & 0 &
\ast & 0 & \ast & 0 & \ast \\  0 & 0 & 0 & 0 & \ast & 0 &
\ast & 0 & \ast \\  z_5 & 0 & nz_5\bar{c}_1+z_5\bar{c}_2 & 0
& \ast & 0 & \ast & 0 & \ast \\  0 & 0 & 0 & 0 & 0 & 0 & \ast & 0
& \ast \\  z_3 & 0 & nz_3\bar{c}_1+z_3\bar{c}_2 & 0 &
{\bold Z}^4/z_3F_2 & 0 & \ast & 0 & \ast \\  0 & 0 & 0 & 0 & 0 & 0
& 0 & 0 & 0 \\  0 & 0 & 0 & 0 & 0 & 0 & 0 & 0 & \ast \\  1
& 0 & Z & 0 & {\bold Z}^4/F_2 & 0 & {\bold Z}^6/F_4 & 0 & {\bold Z}^{10}/F_6
  \endmatrix
$$
We have
$$
d_4(z_3)=\rho^*(\bar{x}_3-\bar{y}_3)=
\rho^*(\sigma_2(\bar{b}_1,\ldots,\bar{b}_5)-
\sigma_2(\bar{a}_1,\ldots,\bar{a}_5))=
$$
$$
=\sigma_2\cdot
\bar{c}_1^2-\sigma_2(\bar{c}_2+\bar{c}_3,\bar{c}_2+
\bar{c}_4,-\bar{c}_3,-\bar{c}_4)=
\sigma_2\cdot
\bar{c}_1^2-\bar{c}_2^2+\bar{c}_3^2+\bar{c}_4^2+
(\bar{c}_3+\bar{c}_4)\bar{c}_2=
$$
$$
=\sigma_2\cdot \bar{c}_1^2-\bar{c}_2^2+\bar{f}.
$$
'Thus, we conclude that
${\bold Z}^4/(F_2\oplus
{\bold Z}(d_4(z_3)))={\bold Z}(\bar{f},
(n-n^2)\bar{c}_1^2+\bar{c}_2^2)/{\bold Z}(d_4(z_3))=
{\bold Z}((n-n^2)\bar{c}_1^2+\bar{c}_2^2)$.
" Now we deduce that
$$
d_4(nz_3\bar{c}_1+z_3\bar{c}_2)=(\sigma_2\cdot
\bar{c}_1^2-\bar{c}_2^2+\bar{f})(n\bar{c}_1+\bar{c}_2)=
$$
$$
=n\sigma_2\cdot \bar{c}_1^3+\sigma_2\cdot \bar{c}_1^2\bar{c}_2-n\cdot
\bar{c}_1\bar{c}_2^2-\bar{c}_2^3+n\cdot
\bar{c}_1\bar{f}+\bar{c}_2\bar{f}.
$$
ŒOne can see that the last element does not vanish in
${\bold Z}^6/F_4$.
Denote by $F_1$ the subgroup, of $H^6B_P$, generated by this
this element. Denote by  $F_2'$ the image of $d_4^{4,3}$
and notice that
$F_2'$ is generated by elements  $d_4(z_3\bar{f})$ and
$d_4((n-n^2)z_3\bar{c}_1^2+z_3\bar{c}_2^2)$. Nondifficult computations show
that  $Ker(d_4^{4,3})=0$.   Let consider $E_5=E_6$.
$$
\matrix z_7 & 0 & \ast & 0 &
\ast & 0 & \ast & 0 & \ast \\  0 & 0 & 0 & 0 & \ast & 0 &
\ast & 0 & \ast \\  z_5 & 0 & nz_5\bar{c}_1+z_5\bar{c}_2 & 0
& \ast & 0 & \ast & 0 & \ast \\  0 & 0 & 0 & 0 & 0 & 0 & \ast & 0
& \ast \\  0 & 0 & 0 & 0 & 0 & 0 & \ast & 0 & \ast \\  0 &
0 & 0 & 0 & 0 & 0 & 0 & 0 & 0 \\  0 & 0 & 0 & 0 & 0 & 0 & 0 &
0 & \ast \\ 1 & 0 & {\bold Z} & 0 & {\bold Z} & 0 &
{\bold Z}^6/(F_4\oplus F_1) & 0 & {\bold Z}^{10}/(F_6\oplus F_2')
\endmatrix
$$
We have
$$
d_6(z_5)=\rho^*(\bar{x}_5-\bar{y}_5)=
\rho^*(\sigma_3(\bar{b}_1,\ldots,\bar{b}_5)-
\sigma_3(\bar{a}_1,\ldots,\bar{a}_5))=
$$
$$
=\sigma_3\cdot
\bar{c}_1^3-\sigma_3(\bar{c}_2+\bar{c}_3,
\bar{c}_2+\bar{c}_4,-\bar{c}_3,-\bar{c}_4)=
$$
$$
=\sigma_3\cdot
\bar{c}_1^3+(\bar{c}_2+\bar{c}_3)(\bar{c}_2+\bar{c}_4)\bar{c}_3+
(\bar{c}_2+\bar{c}_3)(\bar{c}_2+\bar{c}_4)\bar{c}_4-
(\bar{c}_2+\bar{c}_3)\bar{c}_3\bar{c}_4-
(\bar{c}_2+\bar{c}_4)\bar{c}_3\bar{c}_4=
$$
$$
=\sigma_3\cdot
\bar{c}_1^3+(\bar{c}_3+\bar{c}_4)\bar{c}_2^2+
(\bar{c}_3^2+\bar{c}_4^2)\bar{c}_2=
\sigma_3\cdot \bar{c}_1^3+\bar{f}\bar{c}_2.
$$
Let the element $d_6(z_5)$ generates the subgroup
$F_1'$ in $H^6B_P$.
" In addition, denote by $F_1''$ the subgroup generated by the element
$d_6(nz_5\bar{c}_1+z_5\bar{c}_2)$. 
We deduce by simple computation which we omit that
subgroups  $F_1'$ and $F_1''$ are nontrivial in
${\bold Z^6}/(F_4\oplus F_1)$ and  ${\bold Z^{10}}/(F_6\oplus F_2')$,
respectively.
Next,   we consider the  term $E_7$.
$$
\matrix z_7 & 0
& \ast & 0 & \ast & 0 & \ast & 0 & \ast \\  0 & 0 & 0 & 0 & \ast
& 0 & \ast & 0 & \ast \\  0 & 0 & 0 & 0 & \ast & 0 & \ast & 0 &
\ast \\  0 & 0 & 0 & 0 & 0 & 0 & \ast & 0 & \ast \\  0 & 0
& 0 & 0 & 0 & 0 & \ast & 0 & \ast \\  0 & 0 & 0 & 0 & 0 & 0 & 0 &
0 & 0 \\  0 & 0 & 0 & 0 & 0 & 0 & 0 & 0 & \ast \\  1 & 0 &
{\bold Z} & 0 & {\bold Z} & 0 & {\bold Z}^6/(F_4\oplus F_1\oplus F_1') & 0 &
{\bold
Z}/(F_6\oplus F_2'\oplus F_1'')  \endmatrix
$$
' Since
$d_8(z_7)=\sigma_4\bar{c}_1^4-\bar{g}$, the element $z_7$ does not
survive in the next dimensions.
Thus we obtain that
$$
H^1=H^3=H^5=H^7=0,H^2=H^4={\bold Z},
$$
$$
H^6(M_{\bar{p}})=\frac{{\bold Z}^6}{F_4\oplus F_1\oplus F_1'}.
$$
 Let now find $r$ which is equal to the order of group
$H^6$:
$$
r={\Big |} det \left(
\matrix
 \sigma_1 & -2 & 0 & 0 & 0 & 0 \\ 0 & \sigma_1 &
-2 & 0 & 0 & 0 \\ 0 & 0 & \sigma_1 & -2 & 0 & 0 \\ 0 & 0 & 0 & 0 &
\sigma_1 & -2 \\ n\sigma_2 & \sigma_2 & -n & -1 & n & 1 \\ \sigma_3 & 0
& 0 & 0 & 0 & 1
\endmatrix \right){\Big |}.
$$
We remind that
$n=(1-\sigma_1)/2$. By Nondifficult computations, which we omit here, we
obtain that
$$
r=|\sigma_1^3-4\sigma_1\sigma_2+8\sigma_3|.
$$

'Theorem  3 is proved.


\newpage


\Refs

\widestnumber\key{AAA}

\ref
\key Be
\by Berger M.
\paper Les vari\'et\'es Riemannienes homog\`enes normales simplement
connexes \`a courbure strictment positive.
\yr 1961
\vol 15
\pages 179--246
\jour Ann. Scuola Norm. Sup. Pisa
\endref


\ref
\key W
\by Wallach N.R.
\paper Compact homogeneous Riemannian manifolds with strictly
positive curvature.
\yr 1972
\vol 96
\pages 277--295
\jour Ann. of Math.
\endref


\ref
\key AW
\by Aloff S. , Wallach N.R.
\paper  An infinite family of distinct 7-manifolds
admitting positively curved Riemannian structures.
\yr 1975
\vol 81
\pages 93--97
\jour Bull. Amer. Math. Soc.
\endref


\ref
\key BB
\by Berard Bergery L.
\paper  Les vari\'etes Riemannienes homog\'enes simplement
connexes de dimension impair \`a courbure strictment positive.
\yr 1976
\vol 55
\pages 47--68
\jour  J. Pure Math. Appl.
\endref


\ref
\key KS
\by Kreck M., Stolz S.
\paper  Some nondiffeomorphic homeomorphic homogeneous
7--manifolds with positive sectional curvature.
\yr 1991
\vol 33
\pages 465--486
\jour J. Diff. Geom.
\endref


\ref
\key E1
\by Eschenburg J.-H.
\paper  New examples of manifolds with strictly positive
curvature.
\yr 1982
\vol 66
\pages 469--480
\jour Invent. Math.
\endref


\ref
\key E2
\by Eschenburg J.H.
\paper Inhomogeneous spaces of positive curvature.
\yr 1992
\vol 2
\pages 123--132
\jour Differential Geometry and its Applications
\endref


\ref
\key ON
\by O'Neill B.
\paper  The fundamental equations of submersion.
\yr 1966
\vol 23
\pages 459--469
\jour Michigan Math. J.
\endref


\ref
\key M
\by Milnor J.
\book Morse Theory.
\yr 1963
\publ  Ann. of Math. Studies, n. 51, Princeton Univ.
Press, Princeton N.J.
\endref

\ref
\key Bo
\by Borel A.
\paper Sur la cohomologie des espaces fibr\'es principaux et des
espaces homog\`enes de groupes de Lie compacts.
\yr 1953
\vol 57
\pages 251--281
\jour  Ann. of Math.
\endref

\endRefs
\enddocument